\documentclass[11pt]{article}

\usepackage{amsmath}
\usepackage{amssymb}
\usepackage{amsthm}
\usepackage{indentfirst}
\usepackage{enumerate}
\usepackage{cite}
\usepackage{bigints}
\usepackage{mathdots}
\usepackage{arydshln}
\usepackage{multirow}
\usepackage{color}

\newtheorem{dfn}{Definition}[section]
\newtheorem{thm}{Theorem}[section]

\usepackage[OT2,T1]{fontenc}
\input{cyracc.def}

\numberwithin{equation}{section}

\title{\textbf{NOVEL INVARIANTS FOR ALMOST GEODESIC MAPPINGS OF THE THIRD TYPE}}
\author{Du\v san J. Simjanovi\'c${}^{(a)}$, Nenad O. Vesi\'c${}^{(b)}$}
\date{}

\makeatletter
\def\maketag@@@#1{\hbox{\m@th\normalfont\normalsize#1}}
\makeatother
\newcommand\blfootnote[1]{%
  \begingroup
  \renewcommand\thefootnote{}\footnote{#1}%
  \addtocounter{footnote}{-1}%
  \endgroup
}
\begin{document}
  \maketitle

\blfootnote{${}^{(a)}$ Metropolitan University, Tadeu\v sa Ko\v
s\'cu\v ska 63, 11158 Belgrade, Faculty of Information Technology,
Bulevar Sv. Cara Konstantina 80A, 18116 Ni\v s, Serbia}
\blfootnote{${}^{(b)}$ Mathematical Institute of Serbian Academy of
Sciences and Arts. This paper is financially supported by Serbian
Minsitry of Education, Science and Technological Development through
the Mathematical Institute of Serbian Academy of Sciences and Arts}

  \begin{abstract}
    Two kinds of invariance for geometrical objects under
transformations are involved in this paper. With respect to these
kinds, we obtained novel invariants for almost geodesic mappings of
the
    third type of a non-symmetric affine connection space in this
    paper. Our results are presented in two sections. In the Section
    3, we obtained the invariants for the equitorsion almost geodesic mappings
    which do not have the property of reciprocity.
    \\[3pt]


    \textbf{$2010$ Math. Subj. Classification:} 53B05, 15A72, 53A55,
    53A25
  \end{abstract}

  \section{Introduction}

  Manu authors have obtained invariants for different mappings
  between symmetric and non-symmetric affine connection spaces.
  Some
  of them are J. Mike\v s and his research group \cite{mikespi1232015, mikpi1232014,mikpi21999, mikespi22020,miknovi2019}, N. S. Sinyukov \cite{sinjukov}, M.
  S. Stankovi\'c \cite{petrovicmicapi1,micapi3,micapi3spec,jaljubicamicapi3}, M. Z. Petrovi\'c \cite{petrovicmicapi1,
  petrovicmiskolc} and many
  others.

  \subsection{Affine connection spaces}

  An $N$-dimensional manifold $\mathcal M_N$ equipped with a non-symmetric affine
  connection $\nabla$ \big(see Eisenhart \cite{eisNRG}\big) is the non-symmetric affine connection space
  $\mathbb{GA}_N$ \big(in Eisenhart's sense\big).

  The affine connection coefficients of the space $\mathbb{GA}_N$
  are $L^i_{jk}$, $L^i_{jk}\not\equiv L^i_{kj}$ for at least one pair
  $(j,k)\in\{1,\ldots,N\}\times\{1,\ldots,N\}$. For this reason, the
  symmetric and anti-symmetric part of the affine connection
  coefficients $L^i_{jk}$ are defined as

  \begin{equation}
    \begin{array}{ccc}
      L^i_{\underline{jk}}=\dfrac12\big(L^i_{jk}+L^i_{kj}\big)&\mbox{and}&
      L^i_{\underset\vee{jk}}=\dfrac12\big(L^i_{jk}-L^i_{kj}\big),
    \end{array}\label{eq:Lsimantisim}
  \end{equation}

  \noindent respectively.

  The manifold $\mathcal M_N$ equipped with the affine connection
  $\overset0\nabla$ whose coefficients are $L^i_{\underline{jk}}$ is
  the associated space $\mathbb A_N$ \big(of the space
  $\mathbb{GA}_N$\big).

  With respect to the affine connection $\overset0\nabla$, one kind
  of covariant derivative exists
  \cite{sinjukov,mikespi1232015,mikpi21999,mikespi22020,miknovi2019,mikpi1232014}

  \begin{equation}
    a^i_{j|k}=a^i_{j,k}+L^i_{\underline{\alpha
    k}}a^\alpha_j-L^\alpha_{\underline{jk}}a^i_\alpha,
    \label{eq:covderivativeAN}
  \end{equation}

  \noindent for the partial derivative $\partial/\partial x^k$
  denoted by comma.

  In this case, it exists one Ricci-Type identity $a^i_{j|mn}-a^i_{j|nm}=a^\alpha_jR^i_{\alpha mn}-a^i_\alpha
  R^\alpha_{jmn}$, for the curvature tensor

  \begin{equation}
    R^i_{jmn}=L^i_{\underline{jm},n}-L^i_{\underline{jn},m}
    +L^\alpha_{\underline{jm}}L^i_{\underline{\alpha n}}
    -L^\alpha_{\underline{jn}}L^i_{\underline{\alpha m}}.
    \label{eq:RAN}
  \end{equation}

  \noindent of the space $\mathbb A_N$.

  The geometrical object

  \begin{equation}
    R_{ij}=R^\alpha_{ij\alpha},
    \label{eq:RiccTensor}
  \end{equation}

  \noindent is the tensor of the Ricci-curvature.

  Based on the non-symmetric affine connection $\nabla$, S.
  M. Min\v ci\'c defined four kinds of covariant derivatives
  \cite{mincic4,mincic2,mincic1}

  \begin{equation}
    \begin{array}{cc}
      a^i_{j\underset1|k}=a^i_{j,k}+L^i_{\alpha
      k}a^\alpha_j-L^\alpha_{jk}a^i_\alpha,&
      a^i_{j\underset2|k}=a^i_{j,k}+L^i_{k\alpha}a^\alpha_j-L^\alpha_{kj}a^i_\alpha,\\
      a^i_{j\underset3|k}=a^i_{j,k}+L^i_{\alpha
      k}a^\alpha_j-L^\alpha_{kj}a^i_\alpha,&
      a^i_{j\underset4|k}=a^i_{j,k}+L^i_{k\alpha}a^\alpha_j-L^\alpha_{jk}a^i_\alpha.
    \end{array}\label{eq:covderivativeGAN}
  \end{equation}

  With respect to the identities of Ricci-type $a^i_{j\underset
  p|m\underset q|n}-a^i_{j\underset r|n\underset s|m}$,
  $p,q,r,s\in\{0,\ldots,4\}$, $a^i_{j\underset0|k}=a^i_{j|k}$, it is
  obtained the family of curvature tensors for the space
  $\mathbb{GA}_N$

  \begin{equation}
    K^i_{jmn}=R^i_{jmn}+uL^i_{\underset\vee{jm}|n}+u'L^i_{\underset\vee{jn}|m}
    +vL^\alpha_{\underset\vee{jm}}L^i_{\underset\vee{\alpha n}}
    +v'L^\alpha_{\underset\vee{jn}}L^i_{\underset\vee{\alpha m}}
    +wL^\alpha_{\underset\vee{mn}}L^i_{\underset\vee{\alpha j}},
    \label{eq:KGAN}
  \end{equation}

  \noindent for the curvature tensor $R^i_{jmn}$ of the associated space $\mathbb A_N$ and the
  coefficients \cite{mincic4,mincic2,mincic1} $u$, $u'$, $v$, $v'$, $w$.

  The corresponding family of the Ricci-curvatures
  $K_{ij}=K^\alpha_{ij\alpha}$ is

  \begin{equation}
    K_{ij}=R_{ij}+uL^\alpha_{\underset\vee{ij}|\alpha}+
    u'L^\alpha_{\underset\vee{i\alpha}|j}+vL^\alpha_{\underset\vee{ij}}L^\beta_{\underset\vee{\alpha\beta}}-
    (v'+w)L^\alpha_{\underset\vee{i\beta}}L^\beta_{\underset\vee{j\alpha}}.
    \label{eq:RicciGAN}
  \end{equation}




  \subsection{Geodesic and almost geodesic lines and almost geodesic mappings}

  A curve $\ell=\ell(t)$ in the associated space $\mathbb{A}_N$ whose tangential vector
  $\lambda=\dfrac{d\ell}{dt}$ satisfies the system of differential
  equations \cite{sinjukov, miknovi2019}

  \begin{equation}
    \dfrac{d\lambda^i}{dt}+L^i_{\underline{\alpha\beta}}{\lambda^\alpha}{\lambda^\beta}=\zeta\lambda^i,
    \label{eq:geodesiclines}
  \end{equation}

  \noindent for a scalar $\zeta$, is the geodesic line of the space
  $\mathbb{A}_N$.

  A curve $\tilde{\overline \ell}=\tilde{\overline\ell}(t)$ in the associated space
  $\overline{\mathbb A}{}_N$ whose tangential vector
  $\tilde{\overline\lambda}=\dfrac{d\tilde{\overline\lambda}}{dt}$ satisfies the system of
  equations \cite{sinjukov,
  miknovi2019, mikespi1232015,mikpi21999, mikespi22020,mikpi1232014}

  \begin{equation}
    \begin{array}{ccc}
      \tilde{\overline\lambda}{}^i_2=\overline
      a(t)\tilde\lambda{}^i+\overline
      b(t)\tilde{\overline\lambda}{}^i_1,&
      \tilde{\overline\lambda}{}^i_1=\tilde\lambda{}^i_{\|\alpha}\tilde\lambda{}^\alpha,&
      \tilde{\overline\lambda}{}^i_2=\tilde\lambda{}^i_{(1)\|\alpha}\tilde\lambda^\alpha,
    \end{array}\label{eq:almostgeodesicAN}
  \end{equation}

  \noindent for the functions $\overline a(t)$, $\overline b(t)$ and
  the covariant derivative with respect to the affine connection of
  the space $\overline{\mathbb A}{}_N$ denoted by $\|$, is the
  almost geodesic line of the space $\mathbb{\overline A}{}_N$.

  A mapping $f:\mathbb A_N\to\mathbb{\overline A}{}_N$ that any
  geodesic line of the space $\mathbb A_N$ transforms to an almost
  geodesic line of the space $\mathbb{\overline A}{}_N$ is the
  almost geodesic mapping of the space $\mathbb A_N$ \big(see \cite{sinjukov, mikespi1232015,mikpi21999, mikespi22020, miknovi2019,mikpi1232014}\big).

  Three types $\pi_1$, $\pi_2$, $\pi_3$ of almost
  geodesic mappings are detected.

  Sinyukov proved that the inverse mapping of the almost geodesic
  one\linebreak $f:\mathbb A_N\to\mathbb{\overline A}{}_N$ of the type $\pi_3$
  is the almost geodesic mapping of the type $\pi_3$. Almost
  geodesic mappings of the third type whose inverse transformations
  are the almost geodesic mappings of the third type have the property of reciprocity
  \big(see \cite{sinjukov}, page 191\big).

  M. S. Stankovi\'c \cite{micapi1,micapi2,micapi3} generalized the
  Sinyukov's work about almost geodesic mappings. The curve
  $\tilde{\overline\ell}=\tilde{\overline\ell}(t)$ whose tangential
  vector
  $\tilde{\overline\lambda}=\dfrac{d\tilde{\overline\ell}}{dt}$
  is the solution of the system
  \cite{micapi1,micapi2,micapi3,micapi3spec,jaljubicamicapi3,petrovicmicapi1,petrovicmiskolc}

  \begin{equation}
    \begin{array}{ccc}
      \underset k{\tilde{\overline\lambda}}{}^i_{(2)}=
      \underset k{\overline a}(t)\tilde\lambda{}^i+
      \underset k{\overline b}(t)\underset
      k{\tilde{\overline\lambda}}{}^i_{(1)},&
      \underset
      k{\tilde{\overline\lambda}}{}^i_{(1)}=\tilde\lambda{}^i_{\underset k\|\alpha}\tilde\lambda{}^\alpha,&
      \underset k{\tilde{\overline\lambda}}{}^i_{(2)}=\underset
      k{\tilde{\overline\lambda}}{}^i_{(1)\underset
      k\|\alpha}\tilde\lambda{}^\alpha,
    \end{array}\label{eq:almostgeodesicGAN}
  \end{equation}

  \noindent $k\in\{1,2\}$ is the almost geodesic line of the $k$-th
  kind of the space $\mathbb{G\overline A}{}_N$.

  The mapping $f:\mathbb{GA}_N\to\mathbb{G\overline A}{}_N$ that any
  geodesic line of the space $\mathbb{GA}_N$ transforms to the
  almost geodesic line of a $k$-th type, $k=1,2$, of the space
  $\mathbb{G\overline A}{}_N$ is the almost geodesic mapping of the
  space $\mathbb{GA}_N$.

  With respect to the Sinyukov's work, M. S. Stankovi\'c
  \cite{micapi1,micapi2,micapi3} determined three types of almost
  geodesic lines and any of these types is divided into two subtypes.
  These subtypes of almost geodesic mappings are
  $\underset1\pi{}_r$, $\underset2\pi{}_r$, $r=1,2,3$.

  The almost geodesic mapping of a subtype $\underset k\pi{}_3$,
  $k=1,2$,
  has the property of reciprocity if its inverse mapping is
  the almost geodesic mapping of the same subtype.

  We will obtain invariants for almost geodesic mappings of the
  third type of the space $\mathbb{GA}_N$ which have or does not
  have the property of reciprocity below.

  \subsection{Two kinds of invariants for transformations}

  Invariants for mappings between symmetric and non-symmetric affine
  connection spaces are such geometric objects whose values and
  forms do not change under the acting of the corresponding mapping.
  If an almost geodesic mapping
 $f:\mathbb{GA}_N\to\mathbb{G\overline A}{}_N$ \big(or
  $f:\mathbb{A}_N\to\mathbb{\overline A}{}_N$\big) does not have the
  property of reciprocity, the invariants
  of these mappings of the common values and forms can not be obtained. For this reason,
  to obtain the invariants for almost geodesic mappings with respect
  to changes of the curvature tensors under almost geodesic
  mappings, authors assume that these mappings satisfy the property
  of reciprocity.

  Two kinds of invariants are important in physics \big(taken from the textbook \DJ. Mu\v sicki, B. Mili\'c, \emph{Mathematical Foundations of Theoretical
  Physics With a Collection of Solved Problems} \cite{physinvss}, page 103\big):

  \begin{enumerate}[-]
  \item \emph{The covariant} is the object whose form
  stays saved but whose value changes under the transformation of coordinates.
  \item \emph{The invariant} is the object whose form changes but
  whose form stays saved under the transformation of coordinates.
  \item \emph{The total invariant} is the object whose value and
  form stay saved under the transformation.
  \end{enumerate}

  Invariants for mappings between affine connection spaces
  which have been studied in differential geometry
  are analogies to the total invariants from
  physics.  Because all almost geodesic mappings do not have the
  property of reciprocity, we may obtain the geometrical objects
  whose values are preserved unlike their forms under these
  mappings. For this reason, we define two kinds of invariants for
  mappings between the affine connection spaces.

  \begin{dfn}
    Let $f:\mathbb{GA}_N\to\mathbb{G\overline A}{}_N$ be a
    transformation
    and let $U^{i_1\ldots i_p}_{j_1\ldots j_q}$ be a geometrical
    object of the type $(p,q)$.

    \begin{enumerate}[-]
      \item If the transformation $f$ preserves the value of the object $U^{i_1\ldots i_p}_{j_1\ldots
      j_q}$  but
      changes its form to $\overline V{}^{i_1\ldots i_p}_{j_1\ldots
      j_q}$, then the invariance for geometrical object $U^{i_1\ldots i_p}_{j_1\ldots
      j_q}$ under the transformation $f$ is \emph{valued}.
      \item If the transformation $f$ preserves both the value and the form of the geometrical object
      $U^{i_1\ldots i_p}_{j_1\ldots j_q}$, then the invariance for the geometrical object $U^{i_1\ldots i_p}_{j_1\ldots
      j_q}$ under the transformation $f$ is \emph{total}.
    \end{enumerate}
  \end{dfn}

  The transformations which are the main subject of the research in
  this paper are almost geodesic mappings of the third type.

  \subsection{Motivation}

  The transformation rules for affine connection coefficients with
  respect to the third type almost geodesic mappings of symmetric
  and non-symmetric affine connection spaces are

  \begin{equation}
    \begin{array}{ccc}
      \overline
      L{}^i_{\underline{jk}}=L^i_{\underline{jk}}+\psi_j\delta^i_k+\psi_k\delta^i_j&
      \mbox{and}&
      \overline L{}^i_{{jk}}=L^i_{{jk}}+\psi_j\delta^i_k+\psi_k\delta^i_j+\xi^i_{jk},
    \end{array}\label{eq:pi3transformationrules}
  \end{equation}

  \noindent for the $1$-form $\psi_i$, the contravariant vector $\varphi^i$ and the
  tensors $\sigma_{jk}$ and $\xi^i_{jk}$ symmetric and anti-symmetric in the
  covariant indices $j$ and $k$, respectively.

  To generalize the Weyl projective tensor as an invariant for the
  third type almost geodesic mapping of the symmetric affine
  connection space $\mathbb A_N$, N. S. Sinyukov involved \big(see
  \cite{sinjukov}, page 193\big) the geometrical object $q_i$ such
  that $q_\alpha\varphi^\alpha=e$, $e=\pm1$. After some
  computations, Sinyukov generalized the Thomas projective parameter
  and the Weyl projective tensor as the invariants for the almost
  geodesic mapping $f$.

  M. S. Stankovi\'c \big(see \cite{micapi3}\big) continued
  Sinyukov's research about invariants for almost geodesic mappings
  of the third type. With respect to the results presented in
  \cite{micapi3}, N. O. Vesi\'c, Lj. S. Velimirovi\'c and M. S.
  Stankovi\'c \cite{jaljubicamicapi3} obtained the family of invariants for
  equitorsion third type almost geodesic mappings of a non-symmetric
  affine connection space.

  N. S. Sinyukov \cite{sinjukov} obtained one generalization of the Weyl
  projective tensor as invariant for the third type almost geodesic
  mapping of a symmetric affine connection space. M. S. Stankovi\'c
  \cite{micapi3} obtained one generalization of the Weyl projective
  tensor as invariant for the third type almost geodesic mapping of
  a non-symmetric affine connection space with respect to the change
  of the curvature tensor of the corresponding associated space. N.
  O. Vesi\'c, Lj. S. Velimirovi\'c, M. S. Stankovi\'c
  \cite{jaljubicamicapi3} obtained one family of invariants for the
  third type almost geodesic mappings of a non-symmetric affine
  connection space which generalizes the Weyl projective tensor.

  In \cite{jainv1}, it is obtained two invariants for mappings of an
  associated space analogue to the Weyl projective tensor
  \big(called the invariants of the Weyl type\big). That motivated
  us to obtain the invariants for almost geodesic mappings of the
  third type of a non-symmetric affine connection space with respect
  to the results obtained in \cite{jainv1}.

  The formulae of invariants for mappings between symmetric and
  non-symmetric affine connection spaces are obtained in
  \cite{jainv1}. We will use these formulae to meet
  the main goals of this paper. These goals are:

  \begin{enumerate}
    \item To obtain the invariants for equitorsion almost geodesic mappings of
    a non-symmetric affine connection space.
    \item To obtain the necessary and sufficient conditions for
    these invariants to be total.
  \end{enumerate}

  \section{Recall to basic invariants}

  In \cite{jainv1}, the invariants of mappings
  $f:\mathbb{GA}_N\to\mathbb{G\overline A}{}_N$ are obtained. If the
  deformation tensor $P^i_{jk}=\overline L{}^i_{jk}-L^i_{jk}$ of the mapping $f$ is

  \begin{equation}
    P^i_{jk}=\overline
    L{}^i_{jk}-L^i_{jk}=\overline\omega{}^i_{jk}-\omega^i_{jk}+\overline\tau{}^i_{jk}-\tau^i_{jk},
    \label{eq:deformtensorf}
  \end{equation}

  \noindent for
  $\omega^i_{jk}=\omega^i_{kj}$,
  $\overline\omega{}^i_{jk}=\overline\omega{}^i_{kj}$,
  $\tau^i_{jk}=-\tau^i_{kj}$,
  $\overline\tau{}^i_{jk}=-\overline\tau{}^i_{kj}$, the basic
  associated invariants of Thomas and the Weyl type for the mapping $f$ are

  \begin{align}
    &\widetilde{\mathcal
    T}{}^i_{jk}=L^i_{\underline{jk}}-\omega^i_{jk},\label{eq:Thomasbasicsim2}\\\displaybreak[0]
    &\widetilde{\mathcal
    W}{}^i_{jmn}=R^i_{jmn}-\omega^i_{jm|n}+\omega^i_{jn|m}+
    \omega^\alpha_{jm}\omega^i_{\alpha n}-
    \omega^\alpha_{jn}\omega^i_{\alpha m}.\label{eq:Weylbasicsim2}
  \end{align}

  To simplify the last formulae, the next geometrical object
  is used \cite{jainv1}

  \begin{equation}
    L^i_{\underline{jm}|n}=L^i_{\underline{jm},n}+L^i_{\underline{\alpha
    n}}L^\alpha_{\underline{jm}}-L^\alpha_{\underline{jn}}L^i_{\underline{\alpha
    m}}-L^\alpha_{\underline{mn}}L^i_{\underline{j\alpha}}.\label{eq:L|}
  \end{equation}



  In the case of

  \begin{equation}
  \omega^i_{jk}=\delta^i_k\rho_j+\delta^i_j\rho_k+\sigma^i_{jk},
  \label{eq:omega=rho+rho+sigma}
  \end{equation}

  \noindent for
  $\sigma^i_{jk}=\sigma^i_{kj}$, the invariants for the mapping $f$
  given by the equations
  (\ref{eq:Thomasbasicsim2}, \ref{eq:Weylbasicsim2}) transform to

  \begin{align}
    &\widetilde
    T{}^i_{jk}=L^i_{\underline{jk}}-\sigma^i_{jk}-\dfrac1{N+1}\Big(\big(
    L^\alpha_{\underline{j\alpha}}-\sigma^\alpha_{{j\alpha}}\big)\delta^i_k+
    \big(L^\alpha_{\underline{k\alpha}}-\sigma^\alpha_{k\alpha}\big)\delta^i_j
    \Big),\tag{\ref{eq:Thomasbasicsim2}'}\label{eq:Thomasbasicsim2'}\\\displaybreak[0]
    &\aligned
    \widetilde{
    W}{}^i_{jmn}&=R^i_{jmn}-\delta^i_{[m}\rho_{j|n]}-\delta^i_j\rho_{[m|n]}-\sigma^i_{j[m|n]}
    -\delta^i_{[m}\rho_j\rho_{n]}+\delta^i_{[m}\rho_\alpha\sigma^\alpha_{jn]}+
    \sigma^\alpha_{j[m}\sigma^i_{\alpha n]},
    \endaligned\label{eq:basicrhoomega}
    \end{align}

    The derived invariant of the Weyl type for the mapping $f$ is
    \cite{jainv1}

    \begin{equation}
      \aligned
    \widetilde
    W{}^i_{jmn}&=R^i_{jmn}+\dfrac1{N+1}\delta^i_j\big(R_{[mn]}+\sigma^\alpha_{\alpha[m|n]}\big)+
    \dfrac
    N{N^2-1}\delta^i_{[m}R_{jn]}+\dfrac1{N^2-1}\delta^i_{[m}R_{n]j}\\&-\sigma^i_{jm|n}+\sigma^i_{jn|m}+
    \sigma^\alpha_{jm}\sigma^i_{\alpha
    n}-\sigma^\alpha_{jn}\sigma^i_{\alpha m}\\
    &
    -\dfrac1{N^2-1}\delta^i_m\Big(
    \sigma^\alpha_{\alpha[j|n]}+(N+1)\big(\sigma^\alpha_{jn|\alpha}-\sigma^\alpha_{j\alpha|n}-
    \sigma^\alpha_{jn}\sigma^\beta_{\alpha\beta}+\sigma^\alpha_{j\beta}\sigma^\beta_{n\alpha}\big)
    \Big)\\&
    +\dfrac1{N^2-1}\delta^i_n\Big(
    \sigma^\alpha_{\alpha[j|m]}+(N+1)\big(\sigma^\alpha_{jm|\alpha}-\sigma^\alpha_{j\alpha|m}-
    \sigma^\alpha_{jm}\sigma^\beta_{\alpha\beta}+\sigma^\alpha_{j\beta}\sigma^\beta_{m\alpha}\big)
    \Big).
    \endaligned\tag{\ref{eq:Weylbasicsim2}'}\label{eq:Weylbasicsim2'}
    \end{equation}

  The basic invariant for the mapping
  $f:\mathbb{GA}_N\to\mathbb{G\overline A}{}_N$ obtained with
  respect to the transformation rule of the anti-symmetric part
  $L^i_{\underset\vee{jk}}$ of the affine connection coefficient
  $L^i_{jk}$ is

  \begin{equation}
    \hat{\mathcal
    T}{}^i_{jk}=L^i_{\underset\vee{jk}}-\tau^i_{jk}.
    \label{eq:Thomasbasicantisim2}
  \end{equation}

  Let be $\omega^i_{(1).jk}=L^i_{\underline{jk}}$,
  $\overline\omega{}^i_{(1).jk}=\overline L^i_{\underline{jk}}$, $\omega^i_{(2).jk}=\omega^i_{jk}$,
  $\overline\omega{}^i_{(2).jk}=\overline\omega{}^i_{jk}$. The next
  equalities hold \cite{jainv1}

  \begin{align}
    &\aligned
    \hat{\overline{\mathcal T}}{}^i_{jm\|n}-\hat{\mathcal
    T}{}^i_{jm|n}&=
    P^i_{\underline{\alpha n}}\hat{\mathcal T}{}^\alpha_{jm}-
    P^\alpha_{\underline{jn}}\hat{\mathcal T}{}^i_{\alpha m}-
    P^\alpha_{\underline{mn}}\hat{\mathcal T}{}^i_{j\alpha},
    \endaligned\label{eq:torsion-torsioncovariant}\\\displaybreak[0]
    &\aligned
    0&=
    \overline L{}^\alpha_{\underset\vee{jm}}\overline
    L{}^i_{\underset\vee{\alpha
    n}}-L^\alpha_{\underset\vee{jm}}L^i_{\underset\vee{\alpha n}}-
    \overline
    L{}^\alpha_{\underset\vee{jm}}\overline\tau{}^i_{\alpha
    n}-\overline L{}^i_{\underset\vee{\alpha
    n}}\overline\tau{}^\alpha_{jm}+\overline\tau{}^\alpha_{jm}\overline\tau{}^i_{\alpha
    n}\\&+
    L{}^\alpha_{\underset\vee{jm}}\tau{}^i_{\alpha
    n}+L{}^i_{\underset\vee{\alpha
    n}}\tau{}^\alpha_{jm}-\tau{}^\alpha_{jm}\tau{}^i_{\alpha
    n}.
    \endaligned\label{eq:TT-TT=0}
  \end{align}

  Because
  $P^i_{\underline{jk}}=\overline\omega{}^i_{(1).jk}-\omega^i_{jk}=
  \overline\omega{}^i_{(2).jk}-\omega^i_{(2).jk}$, and with respect
  to the equations (\ref{eq:torsion-torsioncovariant},
  \ref{eq:TT-TT=0}), it is obtained \cite{jainv1}

  \begin{equation*}
    \begin{array}{ccc}
      \overline\theta{}^i_{(p).jmn}=\theta^i_{(p).jmn}&\mbox{and}&
      \overline\Theta{}^i_{jmn}=\Theta^i_{jmn},
    \end{array}
  \end{equation*}

  \noindent where $p=(p_1,p_2,p_3)$, $p_1,p_2,p_3\in\{1,2\}$, where

  \begin{align}
    &\theta{}^i_{(p).jmn}=L^i_{\underset\vee{jm}|n}-\tau^i_{jm|n}-
    \omega^i_{(p_1).\alpha n}\hat{\mathcal T}{}^\alpha_{jm}+
    \omega^\alpha_{(p_2).jn}\hat{\mathcal T}{}^i_{\alpha m}+
    \omega^\alpha_{(p_3).mn}\hat{\mathcal
    T}{}^i_{j\alpha},\label{eq:torsioncovariantinv}\\\displaybreak[0]
    &\Theta^i_{jmn}=L^\alpha_{\underset\vee{jm}}L^i_{\underset\vee{\alpha
    n}}-L^\alpha_{\underset\vee{jm}}\tau^i_{\alpha
    n}-L^i_{\underset\vee{\alpha
    n}}\tau^\alpha_{jm}+\tau^\alpha_{jm}\tau^i_{\alpha
    n},\label{eq:torsiontorsioninv}
  \end{align}

  \noindent for the corresponding $\overline\theta{}^i_{(p).jmn}$,
  $\overline\Theta{}^i_{jmn}$.

  The family of invariants for the mapping $f$ with respect to the
  transformation of the family $K^i_{jmn}$ of the curvature tensors
  for the space $\mathbb{GA}_N$ is \cite{jainv1}

  \begin{equation}
    \aligned
    W{}^i_{(p^1).(p^2).jmn}=\widetilde{\mathcal
    W}{}^i_{jmn}+u\theta^i_{(p^1).jmn}+u'\theta^i_{(p^2).jmn}+
    v\Theta^i_{jmn}+v'\Theta^i_{jnm}+w\Theta^i_{mnj},\label{eq:invgeneral}
    \endaligned
  \end{equation}

  \noindent for $p^1=(p^1_1,p^1_2,p^1_3)$,
  $p^2=(p^2_1,p^2_2,p^2_3)$, $p^i_j\in\{1,2\}$ and the corresponding
  invariants $\theta^i_{(p^1).jmn}$, $\theta^i_{(p^2).jmn}$,
  $\Theta^i_{jmn}$ given by (\ref{eq:torsioncovariantinv},
  \ref{eq:torsiontorsioninv}).

  If the mapping $f:\mathbb{GA}_N\to\mathbb{G\overline A}{}_N$ is equitorsion, the invariant $\mathcal{\hat
  T}{}^i_{jk}$ given by (\ref{eq:Thomasbasicantisim2}) reduces to

  \begin{equation}
    \hat{\mathcal T}{}^i_{0.jk}=L^i_{\underset\vee{jk}}.
    \tag{\ref{eq:Thomasbasicantisim2}'}\label{eq:Thomasbasicantisim2'}
  \end{equation}

  The family of invariants of the Weyl type for the equitorsion mapping
  $f:\mathbb{GA}_N\to\mathbb{G\overline A}{}_N$ is

  \begin{equation}
    \mathcal W^i_{0.(p^1).(p^2).jmn}=\widetilde{\mathcal
    W}{}^i_{jmn}+u\theta^i_{0.(p^1).jmn}+u'\theta^i_{0.(p^2).jmn},\label{eq:equitorsioninvgeneral}
  \end{equation}

  \noindent where

  \begin{equation}
    \theta^i_{0.(p).jmn}=L^i_{\underset\vee{jm}|n}-\omega^i_{(p_1).\alpha
    n}L^\alpha_{\underset\vee{jm}}+\omega^\alpha_{(p_2).jn}L^i_{\underset\vee{\alpha
    m}}+\omega^\alpha_{(p_3).mn}L^i_{\underset\vee{j\alpha}}.
    \tag{\ref{eq:torsioncovariantinv}'}\label{eq:torsioncovariantinv'}
  \end{equation}

  \section{Invariants for equitorsion almost geodesic mappings}

  Let $f:\mathbb{GA}_N\to\mathbb{G\overline A}{}_N$ be an equitorsion almost
  geodesic mapping of the third type. Its basic equations are
  \cite{micapi3}

  \begin{equation}
    \left\{
    \begin{array}{l}
      \overline
      L{}^i_{\underline{jk}}=L^i_{\underline{jk}}+\psi_j\delta^i_k+\psi_k\delta^i_j+2\sigma_{jk}\varphi^i,\\
      \varphi^i_{\underset1{|}j}=\underset1{\nu}{}_j\varphi^i+\underset1\mu\delta^i_j.
    \end{array}
    \right.\label{eq:pi31basic}
  \end{equation}

  Let us rewrite the first of the last basic equations as

  \begin{equation}
    \overline
    L{}^i_{\underline{jk}}=L^i_{\underline{jk}}+\psi_j\delta^i_k+\psi_k\delta^i_j+D^i_{jk},
    \tag{\ref{eq:pi31basic}'}\label{eq:pi31basic1'}
  \end{equation}

  \noindent for the tensor $D^i_{jk}$,
  $D^i_{jk}=D^i_{kj}=2\sigma_{jk}\varphi^i$.

  In the case of the inverse mapping $f^{-1}:\mathbb{G\overline
  A}{}_N\to\mathbb{GA}_N$, it exists the tensor $\overline
  D{}^i_{jk}$, $\overline D{}^i_{jk}=\overline D{}^i_{kj}$, $\overline D{}^i_{jk}=-D^i_{jk}$, such
  that

  \begin{equation}
    \aligned
    L^i_{\underline{jk}}=\overline
    L{}^i_{\underline{jk}}-\psi_j\delta^i_k-\psi_k\delta^i_j-D^i_{jk}
    =\overline L{}^i_{\underline{jk}}-\psi{}_j\delta^i_k-\psi{}_k\delta^i_j+\overline
    D{}^i_{jk}.
    \endaligned\label{eq:LtoLinversepi31}
  \end{equation}

  Hence, the equation (\ref{eq:pi31basic1'})
  transforms to

  \begin{equation}
    \overline
    L{}^i_{\underline{jk}}=L^i_{\underline{jk}}+\psi_j\delta^i_k+\psi_k\delta^i_j-\dfrac12\big(\overline D{}^i_{jk}-
    D^i_{jk}\big).
    \tag{\ref{eq:pi31basic}''}\label{eq:pi31basic'a}
  \end{equation}

  After contracting the last equation by $i$ and $k$, one gets

  \begin{equation}
    \psi_j=\dfrac1{N+1}\big(\overline
    L{}^\alpha_{\underline{j\alpha}}+\dfrac12\overline D{}^\alpha_{j\alpha}\big)
    -\dfrac1{N+1}\big(L^\alpha_{\underline{j\alpha}}+\dfrac12D^\alpha_{j\alpha}\big).
    \label{eq:psijpi31}
  \end{equation}

  If substitute the expression (\ref{eq:psijpi31}) in the
  equation (\ref{eq:pi31basic'a}) and use the expression $D^i_{jk}=2\sigma_{jk}\varphi^i$, we will obtain

  \begin{equation}
    \omega^i_{jk}=\dfrac1{N+1}\delta^i_k\big(L^\alpha_{\underline{j\alpha}}+\sigma_{j\alpha}\varphi^\alpha\big)+
    \dfrac1{N+1}\delta^i_j\big(L^\alpha_{\underline{k\alpha}}+\sigma_{k\alpha}\varphi^\alpha\big)-\sigma_{jk}\varphi^i.
    \label{eq:pi31omega}
  \end{equation}

  The second of the basic equations (\ref{eq:pi31basic}) is
  equivalent to

  \begin{equation}
    \varphi^i_{|j}=-L^i_{\underset\vee{\alpha
    j}}\varphi^\alpha+\underset1\nu{}_{j}\varphi^i+\underset1\mu\delta^i_j.
    \label{eq:pi31basic'b}
  \end{equation}

  After comparing the equations (\ref{eq:omega=rho+rho+sigma}) and
  (\ref{eq:pi31basic'b}), one reads

  \begin{equation}
  \begin{array}{ccc}\rho_j=\dfrac1{N+1}\big(L^\alpha_{\underline{j\alpha}}+\sigma_{j\alpha}\varphi^\alpha\big)&
  \mbox{and}&
  \sigma^i_{jk}=-\sigma_{jk}\varphi^i.
  \end{array}\label{eq:rho=sigma=pi31}
  \end{equation}

  Hence, we get

  \begin{equation}
    \left\{\begin{array}{l}
    -\sigma^i_{jm|n}=\big(\sigma_{jm}\varphi^i\big)_{|n}=\sigma_{jm|n}\varphi^i-\sigma_{jm}L^i_{\underset\vee{\alpha n}}\varphi^\alpha
    +\sigma_{jm}\underset1\nu{}_n\varphi^i+\sigma_{jm}\underset1\mu\delta^i_n,\\
    -\sigma^\alpha_{ij|\alpha}=\big(\sigma_{ij}\varphi^\alpha)_{|\alpha}=\sigma_{ij|\alpha}\varphi^\alpha-
    \sigma_{ij}L^\beta_{\underset\vee{\alpha\beta}}\varphi^\alpha+
    \sigma_{ij}\underset1\nu{}_\alpha\varphi^\alpha+N\underset1\mu\sigma_{ij},\\
    -\sigma^\alpha_{\alpha i|j}=\big(\sigma_{\alpha i}\varphi^\alpha\big)_{|j}=
    \big(\sigma_{i\alpha}\varphi^\alpha\big)_{|j}=\sigma_{\alpha
    i|j}\varphi^\alpha-\sigma_{\beta i}L^\beta_{\underset\vee{\alpha
    j}}\varphi^\alpha+\sigma_{\alpha
    i}\underset1\nu{}_j\varphi^\alpha+
    \underset1\mu\sigma_{ij}.
    \end{array}\right.\label{eq:sigmaijm|n1}
  \end{equation}

  After substituting the expressions (\ref{eq:pi31omega},
  \ref{eq:rho=sigma=pi31}, \ref{eq:sigmaijm|n1}) in the equations
  (\ref{eq:Thomasbasicsim2'}, \ref{eq:Weylbasicsim2'},
  \ref{eq:torsioncovariantinv'}), we obtain the next geometrical
  objects

  \begin{align}
    &\aligned
    \underset1{\widetilde{\mathcal
    T}}{}^i_{jk}=L^i_{\underline{jk}}-\dfrac1{N+1}\delta^i_k\big(L^\alpha_{\underline{j\alpha}}+\sigma_{j\alpha}\varphi^\alpha\big)
    -\dfrac1{N+1}\delta^i_j\big(L^\alpha_{\underline{k\alpha}}+\sigma_{k\alpha}\varphi^\alpha\big)
    +\sigma_{jk}\varphi^i,
    \endaligned\label{eq:basicThomaspi31}\\\displaybreak[0]
    &\aligned
    \underset1{\mathcal{\widetilde
    W}}{}^i_{jmn}&=R^i_{jmn}+\dfrac1{N+1}\delta^i_j\Big(R_{[mn]}-
    \big(\sigma_{[m\alpha|n]}-\sigma_{[m\beta}L^\beta_{\underset\vee{\alpha
    n}]}+\sigma_{[m\alpha}\underset1\nu{}_{n]}\big)\varphi^\alpha\Big)\\&
    +\sigma_{j[m|n]}\varphi^i+\sigma_{j[m}\sigma_{\alpha n]}\varphi^\alpha\varphi^i
    -\sigma_{j[m}L^i_{\underset\vee{\alpha n}]}\varphi^\alpha+
    \sigma_{j[m}\underset1\nu{}_{n]}\varphi^i-\delta^i_{[m}\underset1\mu\sigma_{jn]}\\&-
    \dfrac1{N+1}\Big(\delta^i_{[m}L^\alpha_{\underline{j\alpha}|n]}+\big(\delta^i_{[m}\sigma_{j\alpha|n]}-
    \delta^i_{[m}\sigma_{j\beta}L^\beta_{\underset\vee{\alpha
    n}]}+\delta^i_{[m}\sigma_{j\alpha}\underset1\nu{}_{n]}\big)\varphi^\alpha+
    \delta^i_{[m}\underset1\mu\sigma_{jn]}\Big)\\&+
    \dfrac1{N+1}\delta^i_{[m}\sigma_{jn]}\big(L^\beta_{\underline{\alpha\beta}}+\sigma_{\alpha\beta}\varphi^\beta\big)\varphi^\alpha
    \\&-\dfrac1{(N+1)^2}\big(L^\alpha_{\underline{j\alpha}}+\sigma_{j\alpha}\varphi^\alpha\big)
    \big(\delta^i_{[m}L^\beta_{\underline{n]\beta}}+\delta^i_{[m}\sigma_{n]\beta}\varphi^\beta\big),
    \endaligned\label{eq:Weylbasicsimpi31}\\\displaybreak[0]
    &\aligned
    \underset1{\widetilde
    W}{}^i_{jmn}&=R^i_{jmn}+\dfrac1{N+1}\delta^i_jR_{[mn]}+\dfrac
    N{N^2-1}\delta^i_{[m}R_{jn]}+\dfrac1{N^2-1}\delta^i_{[m}R_{n]j}
    \\&+\sigma_{j[m|n]}\varphi^i-\sigma_{j[m}L^i_{\underset\vee{\alpha
    n}]}\varphi^\alpha+\sigma_{j[m}\underset1\nu{}_{n]}\varphi^i+
    \sigma_{j[m}\sigma_{\alpha
    n]}\varphi^\alpha\varphi^i\\&-\dfrac1{N+1}\delta^i_j\big(\sigma_{\alpha[m|n]}-\sigma_{\beta[m}L^\beta_{\underset\vee{\alpha
    n}]}+\sigma_{\alpha[m}\underset1\nu{}_{n]}\big)\varphi^\alpha\\&+
    \dfrac1{N-1}\big(\delta^i_{[m}\sigma_{jn]|\alpha}-\delta^i_{[m}\sigma_{jn]}L^\beta_{\underset\vee{\alpha\beta}}+
    \delta^i_{[m}\sigma_{jn]}\underset1\nu{}_\alpha\big)\varphi^\alpha\\&+
    \dfrac1{N-1}\big(\delta^i_{[m}\sigma_{jn]}\sigma_{\alpha\beta}-\delta^i_{[m}\sigma_{j\alpha}\sigma_{n]\beta}\big)
    \varphi^\alpha\varphi^\beta\\&
    -\dfrac N{N^2-1}
    \big(\delta^i_{[m}\sigma_{\alpha j|n]}-\delta^i_{[m}\sigma_{\beta j}L^\beta_{\underset\vee{\alpha
    n}]}+\delta^i_{[m}\sigma_{\alpha j}\underset1\nu{}_{n]}\big)\varphi^\alpha\\&
    -\dfrac1{N^2-1}
    \big(\delta^i_{[m}\sigma_{\alpha n]|j}-\delta^i_{[m}\sigma_{\beta n]}L^\beta_{\underset\vee{\alpha
    j}}+\delta^i_{[m}\sigma_{\alpha n]}\underset1\nu{}_{j}\big)\varphi^\alpha.
    \endaligned\label{eq:basicWeylpi31}
  \end{align}

  Let us express the invariant $\underset1{\widetilde W}{}^i_{jmn}$ in the form

  \begin{equation}
    \aligned
    \underset1{\widetilde
    W}{}^i_{jmn}&=R^i_{jmn}+\dfrac1{N+1}\delta^i_jR{}_{[mn]}+
    \dfrac
    N{N^2-1}\delta^i_{[m}R_{jn]}+\dfrac1{N^2-1}\delta^i_{[m}R_{n]j}\\&+
    \delta^i_j\underset1X{}_{[mn]}+\delta^i_{[m}\underset1Y{}_{jn]}+\underset1Z{}^i_{jmn},
    \endaligned\tag{\ref{eq:basicWeylpi31}'}\label{eq:basicWeylpi31'}
  \end{equation}

  \noindent for the corresponding tensors

  \begin{align}
    &\aligned
    \underset1X{}_{ij}&=-\dfrac1{N+1}\big(\sigma_{\alpha i|j}-\sigma_{\beta i}L^\beta_{\underset\vee{\alpha
    j}}+\sigma_{\alpha i}\underset1\nu{}_j\big)\varphi^\alpha,
    \endaligned\label{eq:pi31Xij}\\
    &\aligned
    \underset1Y{}_{ij}&=\dfrac1{N-1}\big(\sigma_{ij|\alpha}-
    \sigma_{ij}L^\beta_{\underset\vee{\alpha\beta}}+
    \sigma_{ij}\underset1\nu{}_\alpha\big)\varphi^\alpha+
    \dfrac1{N-1}\big(\sigma_{ij}\sigma_{\alpha\beta}-\sigma_{i\alpha}\sigma_{j\beta}\big)\varphi^\alpha\varphi^\beta\\&-
    \dfrac N{N^2-1}\big(\sigma_{\alpha i|j}-\sigma_{\beta i}L^\beta_{\underset\vee{\alpha
    j}}+\sigma_{\alpha i}\underset1\nu{}_{j}\big)\varphi^\alpha-
    \dfrac1{N^2-1}\big(\sigma_{\alpha j|i}-\sigma_{\beta j}L^\beta_{\underset\vee{\alpha i}}+
    \sigma_{\alpha j}\underset1\nu{}_i\big)\varphi^\alpha,
    \endaligned\label{eq:pi31Yij}\\
    &\aligned
    \underset1Z{}^i_{jmn}&=\sigma_{j[m|n]}\varphi^i-\sigma_{j[m}L^i_{\underset\vee{\alpha
    n}]}\varphi^\alpha+\sigma_{j[m}\underset1\nu{}_{n]}\varphi^i+
    \sigma_{j[m}\sigma_{\alpha
    n]}\varphi^\alpha\varphi^i.
    \endaligned\label{eq:pi31Zijmn}
  \end{align}

  After contracting the equality $0=\underset1{\widetilde{\overline
  W}}{}^i_{jmn}-\underset1{\widetilde{W}}{}^i_{jmn}$ by the indices $i$ and $j$,
  one gets

  \begin{equation}
    \aligned
    \underset1{\overline X}{}_{[mn]}-\underset1X{}_{[mn]}&=-\dfrac1N\big(\underset1{\overline
    Y}{}_{[mn]}+\underset1{\overline Z}{}^\alpha_{\alpha
    mn}\big)+\dfrac1N\big(\underset1Y{}_{[mn]}+\underset1Z{}^\alpha_{\alpha
    mn}\big),
    \endaligned\label{eq:X-Xpi31}
  \end{equation}

  \noindent where

  \begin{equation}
    \underset1Y{}_{[ij]}=-\dfrac1{N+1}\big(\sigma_{\alpha[i|j]}-\sigma_{\beta[i}L^\beta_{\underset\vee{\alpha
    j}]}+\sigma_{\alpha[i}\underset1\nu{}_{j]}\big)\varphi^\alpha.
    \label{eq:Y[mn]pi31}
  \end{equation}

  With respect to the equations (\ref{eq:pi31Yij},
  \ref{eq:pi31Zijmn}, \ref{eq:X-Xpi31}, \ref{eq:Y[mn]pi31})
  substituted into the equality

  \begin{equation}
    \aligned
    0&=\underset1{\widetilde{\overline W}}{}^i_{jmn}-
    \underset1{\widetilde W}{}^i_{jmn}\\&=\overline
    R{}^i_{jmn}-R^i_{jmn}+\dfrac1{N+1}\delta^i_j\big(\overline
    R{}_{[mn]}-R_{[mn]}\big)+
    \dfrac N{N^2-1}\big(\delta^i_{[m}\overline
    R{}_{jn]}-\delta^i_{[m}R_{jn]}\big)\\&+
    \dfrac1{N^2-1}\big(\delta^i_{[m}\overline R{}_{n]j}-
    \delta^i_{[m}R_{n]j}\big)+
    \delta^i_j\big(\underset1{\overline X}{}_{[mn]}-\underset1X{}_{[mn]}\big)\\&+
    \big(\delta^i_{[m}\underset1{\overline Y}{}_{jn]}-\delta^i_{[m}\underset1Y{}_{jn]}\big)+
    \underset1{\overline Z}{}^i_{jmn}-\underset1Z{}^i_{jmn},
    \endaligned\label{eq:derinv-derinvpi31}
  \end{equation}

  \noindent one gets $\underset1{\widetilde{\widetilde{\overline
  W}}}{}^i_{jmn}=\underset1{\widetilde{\widetilde W}}{}^i_{jmn}$,
  where

  \begin{equation}
    \aligned
    \underset1{\widetilde{\widetilde
    W}}{}^i_{jmn}&=R^i_{jmn}+\dfrac1{N+1}\delta^i_jR_{[mn]}+\dfrac
    N{N^2-1}\delta^i_{[m}R_{jn]}+\dfrac1{N^2-1}\delta^i_{[m}R_{n]j}\\&+
    \sigma_{j[m|n]}\varphi^i-\sigma_{j[m}L^i_{\underset\vee{\alpha
    n}]}\varphi^\alpha+\sigma_{j[m}\underset1\nu{}_{n]}\varphi^i+\sigma_{j[m}\sigma_{\alpha
    n]}\varphi^\alpha\varphi^i\\&-
    \dfrac
    N{N+1}\delta^i_j\big(\sigma_{\alpha[m|n]}-\sigma_{\beta[m}L^\beta_{\underset\vee{\alpha
    n}]}+\sigma_{\alpha[m}\underset1\nu{}_{n]}\big)\varphi^\alpha\\&+
    \dfrac1{N-1}\big(\delta^i_{[m}\sigma_{jn]|\alpha}-
    \delta^i_{[m}\sigma_{jn]}L^\beta_{\underset\vee{\alpha \beta}}+
    \delta^i_{[m}\sigma_{jn]}\underset1\nu{}_\alpha\big)\varphi^\alpha\\&+
    \dfrac1{N-1}\big(\delta^i_{[m}\sigma_{jn]}\sigma_{\alpha\beta}-
    \delta^i_{[m}\sigma_{j\alpha}\sigma_{n]\beta}\big)\varphi^\alpha\varphi^\beta\\&-
    \dfrac N{N^2-1}\big(\delta^i_{[m}\sigma_{\alpha
    j|n]}-\delta^i_{[m}\sigma_{\beta j}L^\beta_{\underset\vee{\alpha
    n}]}+\delta^i_{[m}\underset1\nu{}_{n]}\big)\varphi^\alpha\\&-
    \dfrac1{N^2-1}\big(\delta^i_{[m}\sigma_{\alpha
    n]|j}-\delta^i_{[m}\sigma_{\beta
    n]}L^\beta_{\underset\vee{\alpha j}}+\delta^i_{[m}\sigma_{\alpha
    n]}\underset1\nu{}_j\big)\varphi^\alpha.
    \endaligned\label{eq:inv3pi31}
  \end{equation}

  Based on the invariants $\underset1{\widetilde W}{}^i_{jmn}$ and
  $\underset1{\widetilde{\widetilde W}}{}^i_{jmn}$, one concludes
  that the geometrical object \linebreak
  $\big(\sigma_{\alpha[i|j]}-\sigma_{\beta[i}L^\beta_{\underset\vee{\alpha
  j}]}+\sigma_{\alpha[i}\underset1\nu{}_{j]}\big)\varphi^\alpha$ is
  an invariant for the mapping $f$. Hence, the invariant
  $\underset1{\widetilde W}{}^i_{jmn}$ given by (\ref{eq:basicWeylpi31}) reduces to

  \begin{equation}
    \aligned
    \underset1{\widetilde{
    W}}{}^i_{jmn}&=R^i_{jmn}+\dfrac1{N+1}\delta^i_jR_{[mn]}+\dfrac
    N{N^2-1}\delta^i_{[m}R_{jn]}+\dfrac1{N^2-1}\delta^i_{[m}R_{n]j}\\&+
    \sigma_{j[m|n]}\varphi^i-\sigma_{j[m}L^i_{\underset\vee{\alpha
    n}]}\varphi^\alpha+\sigma_{j[m}\underset1\nu{}_{n]}\varphi^i+\sigma_{j[m}\sigma_{\alpha
    n]}\varphi^\alpha\varphi^i\\&+
    \dfrac1{N-1}\big(\delta^i_{[m}\sigma_{jn]|\alpha}-
    \delta^i_{[m}\sigma_{jn]}L^\beta_{\underset\vee{\alpha \beta}}+
    \delta^i_{[m}\sigma_{jn]}\underset1\nu{}_\alpha\big)\varphi^\alpha\\&+
    \dfrac1{N-1}\big(\delta^i_{[m}\sigma_{jn]}\sigma_{\alpha\beta}-
    \delta^i_{[m}\sigma_{j\alpha}\sigma_{n]\beta}\big)\varphi^\alpha\varphi^\beta\\&-
    \dfrac N{N^2-1}\big(\delta^i_{[m}\sigma_{\alpha
    j|n]}-\delta^i_{[m}\sigma_{\beta j}L^\beta_{\underset\vee{\alpha
    n}]}+\delta^i_{[m}\underset1\nu{}_{n]}\big)\varphi^\alpha\\&-
    \dfrac1{N^2-1}\big(\delta^i_{[m}\sigma_{\alpha
    n]|j}-\delta^i_{[m}\sigma_{\beta
    n]}L^\beta_{\underset\vee{\alpha j}}+\delta^i_{[m}\sigma_{\alpha
    n]}\underset1\nu{}_j\big)\varphi^\alpha.
    \endaligned\label{eq:basicWeylpi31'}
  \end{equation}

  If contracts the equality $0=\underset1{\widetilde{\overline
  W}}{}^i_{jmn}-\underset1{\widetilde{W}}{}^i_{jmn}$ \big(equivalent to the equation
  (\ref{eq:derinv-derinvpi31})\big)
   by the indices $i$ and $n$
  and anti-symmetrizes the contracted equation by the indices $j$ and $m$,
  one will obtain

  \begin{equation}
    \underset1{\overline X}{}_{[jm]}-\underset1X{}_{[jm]}=-\dfrac{N-1}2\big(\underset1{\overline
    Y}{}_{[jm]}-\underset1{Y}{}_{[jm]}\big)+\dfrac12\big(\underset1{\overline
    Z}{}^\alpha_{[jm]\alpha}-\underset1Z{}^\alpha_{[jm]\alpha}\big).\label{eq:X-Xpi31i=n}
  \end{equation}

  After substituting the expression (\ref{eq:X-Xpi31i=n}) in the
  equation (\ref{eq:derinv-derinvpi31}), we obtain
  $\underset1{\widetilde{\widetilde{\widetilde{\overline
  W}}}}{}^i_{jmn}=\underset1{\widetilde{\widetilde{\widetilde W}}}{}^i_{jmn}$, for

  \begin{equation}
    \aligned
    \underset1{\widetilde{\widetilde{\widetilde W}}}{}^i_{jmn}&=
    R^i_{jmn}+\dfrac1{N+1}\delta^i_jR_{[mn]}+\dfrac
    N{N^2-1}\delta^i_{[m}R_{jn]}+\dfrac1{N^2-1}\delta^i_{[m}R_{n]j}\\&+
    \sigma_{j[m|n]}\varphi^i-\sigma_{j[m}L^i_{\underset\vee{\alpha
    n}]}\varphi^\alpha+\sigma_{j[m}\underset1\nu{}_{n]}\varphi^i+\sigma_{j[m}\sigma_{\alpha
    n]}\varphi^\alpha\varphi^i\\&+
    \dfrac{N-1}{2(N+1)}\delta^i_j\big(\sigma_{\alpha[m|n]}-
    \sigma_{\beta[m}L^\beta_{\underset\vee{\alpha
    n}]}+\sigma_{\alpha[m}\underset1\nu{}_{n]}\big)\varphi^\alpha\\&-
    \dfrac12\delta^i_j\big(\sigma_{[m\alpha|n]}-\sigma_{[m\beta}L^\beta_{\underset\vee{\alpha n}]}+\sigma_{[m\alpha}\underset1\nu{}_{n]}\big)\varphi^\alpha\\&+
    \dfrac1{N-1}\big(\delta^i_{[m}\sigma_{jn]|\alpha}-\delta^i_{[m}\sigma_{jn]}L^\beta_{\underset\vee{\alpha\beta}}+
    \delta^i_{[m}\sigma_{jn]}\underset1\nu{}_\alpha\big)\varphi^\alpha\\&+
    \dfrac1{N-1}\big(\delta^i_{[m}\sigma_{jn]}\sigma_{\alpha\beta}-\delta^i_{[m}\sigma_{j\alpha}\sigma_{n]\beta}\big)
    \varphi^\alpha\varphi^\beta\\&
    -\dfrac N{N^2-1}
    \big(\delta^i_{[m}\sigma_{\alpha j|n]}-\delta^i_{[m}\sigma_{\beta j}L^\beta_{\underset\vee{\alpha
    n}]}+\delta^i_{[m}\sigma_{\alpha j}\underset1\nu{}_{n]}\big)\varphi^\alpha\\&
    -\dfrac1{N^2-1}
    \big(\delta^i_{[m}\sigma_{\alpha n]|j}-\delta^i_{[m}\sigma_{\beta n]}L^\beta_{\underset\vee{\alpha
    j}}+\delta^i_{[m}\sigma_{\alpha
    n]}\underset1\nu{}_{j}\big)\varphi^\alpha.
    \endaligned\label{eq:inv4pi31}
  \end{equation}

  After comparing the invariants
  $\underset1{\widetilde{\widetilde{\widetilde W}}}{}^i_{jmn}$ and
  $\underset1{\widetilde W}{}^i_{jmn}$ for the mapping $f$, one gets
  that the invariant $\underset1{\widetilde{\widetilde{\widetilde
  W}}}{}^i_{jmn}$ given by (\ref{eq:inv4pi31}) reduces to the
  invariant $\underset1{\widetilde W}{}^i_{jmn}$ given by the
  equation (\ref{eq:basicWeylpi31'}).

  With respect to the transformation of the family (\ref{eq:KGAN})
  of the curvature tensors of the space $\mathbb{GA}_N$ under the mapping
  $f$, and with respect to the equation
  (\ref{eq:equitorsioninvgeneral}), we obtain the next geometrical
  objects.

  \begin{align}
    &\aligned
    \underset1{\mathcal W}{}^i_{0.(p^1).(p^2).jmn}&=\underset1{\widetilde{\mathcal W}}{}^i_{jmn}
    +uL^i_{\underset\vee{jm}|n}+u'L^i_{\underset\vee{jn}|m}\\&-u\big(\omega^i_{(p_1^1).\alpha n}L^\alpha_{\underset\vee{jm}}-
    \omega^\alpha_{(p_2^1).jn}L^i_{\underset\vee{\alpha m}}-
    \omega^\alpha_{(p_3^1).mn}L^i_{\underset\vee{j\alpha}}\big)\\&-u'\big(
    \omega^i_{(p_1^2).\alpha m}L^\alpha_{\underset\vee{jn}}-
    \omega^\alpha_{(p_2^2).jm}L^i_{\underset\vee{\alpha n}}-
    \omega^\alpha_{(p_3^2).mn}L^i_{\underset\vee{j\alpha}}\big),
    \endaligned\label{eq:inv1KtoKpi31}\\\displaybreak[0]
    &\aligned
    \underset1W{}^i_{0.(p^1).(p^2).jmn}&=\underset1{\widetilde{W}}{}^i_{jmn}
    +uL^i_{\underset\vee{jm}|n}+u'L^i_{\underset\vee{jn}|m}\\&-u\big(\omega^i_{(p_1^1).\alpha n}L^\alpha_{\underset\vee{jm}}-
    \omega^\alpha_{(p_2^1).jn}L^i_{\underset\vee{\alpha m}}-
    \omega^\alpha_{(p_3^1).mn}L^i_{\underset\vee{j\alpha}}\big)\\&-u'\big(
    \omega^i_{(p_1^2).\alpha m}L^\alpha_{\underset\vee{jn}}-
    \omega^\alpha_{(p_2^2).jm}L^i_{\underset\vee{\alpha n}}-
    \omega^\alpha_{(p_3^2).mn}L^i_{\underset\vee{j\alpha}}\big),
    \endaligned\label{eq:inv2KtoKpi31}
  \end{align}

  \noindent for $p^1_1,\ldots,p^2_3\in\{1,2\}$,
  $\omega^i_{(1).jk}=L^i_{\underline{jk}}$ and
  $\omega^i_{(2).jk}=\omega^i_{jk}$, for the geometrical object $\omega^i_{jk}$
  given in the equation (\ref{eq:pi31omega}).

  It holds the next theorem.

  \begin{thm}
    Let $f:\mathbb{GA}_N\to\mathbb{G\overline A}{}_N$ be an equitorsion almost
    geodesic mapping of the type $\underset1\pi{}_3$.

    The geometrical object $\widetilde{\mathcal T}{}^i_{jk}$ given
    by \emph{(\ref{eq:basicThomaspi31})} is the basic invariant of
    the Thomas type for the mapping $f$. The invariance of this
    geometrical object is total.

    The geometrical object $\widetilde{\mathcal W}{}^i_{jmn}$ given
    by \emph{(\ref{eq:Weylbasicsimpi31})} is the basic associated
    invariant for the mapping $f$. The invariance of this
    geometrical object is valued. It is total if and only if the
    mapping $f$ has the property of reciprocity.

    The geometrical object $\widetilde W{}^i_{jmn}$ given by
    \emph{(\ref{eq:basicWeylpi31'})} is the associated derived invariant of the
    Weyl type for the mapping $f$. The invariance of this
    geometrical object is valued. It is total if and only if the
    mapping $f$ has the property of reciprocity.

    The geometrical objects $\underset1{\mathcal W}{}^i_{jmn}$,
    $\underset1W{}^i_{jmn}$, given by
    \emph{(\ref{eq:inv1KtoKpi31}, \ref{eq:inv2KtoKpi31})}, are the invariants
    for the equitorsion third type almost geodesic mapping $f$. The invariance of
    these
    geometrical objects are valued. They are total if and only if the
    mapping $f$ has the property of reciprocity.
    \qed
  \end{thm}

  The basic equations for the almost geodesic mapping
  $f:\mathbb{GA}_N\to\mathbb{G\overline A}{}_N$ of the type
  $\underset2\pi{}_3$ are

  \begin{equation}
    \left\{
    \begin{array}{l}
      \overline
      L{}^i_{\underline{jk}}=L^i_{\underline{jk}}+\psi_j\delta^i_k+\psi_k\delta^i_j+2\sigma_{jk}\varphi^i,\\
      \varphi^i_{\underset2{|}j}=\underset2{\nu}{}_j\varphi^i+\underset2\mu\delta^i_j.
    \end{array}
    \right.\label{eq:pit2basic}
  \end{equation}

  The second of the basic equations (\ref{eq:pit2basic}) is
  equivalent to

  \begin{equation}
    \varphi^i_{|j}=L^i_{\underset\vee{\alpha j}}\varphi^\alpha+
    \underset2{\nu}{}_j\varphi^i+\underset2\mu\delta^i_j.\label{eq:pi32basic'b}
  \end{equation}

  Analogously as above, one obtains the following geometrical
  objects.

  \begin{align}
    &\aligned
    \underset2{\widetilde{\mathcal
    T}}{}^i_{jk}=L^i_{\underline{jk}}-\dfrac1{N+1}\delta^i_k\big(L^\alpha_{\underline{j\alpha}}+\sigma_{j\alpha}\varphi^\alpha\big)
    -\dfrac1{N+1}\delta^i_j\big(L^\alpha_{\underline{k\alpha}}+\sigma_{k\alpha}\varphi^\alpha\big)
    +\sigma_{jk}\varphi^i,
    \endaligned\label{eq:basicThomaspi32}\\\displaybreak[0]
    &\aligned
    \underset2{\mathcal{\widetilde
    W}}{}^i_{jmn}&=R^i_{jmn}+\dfrac1{N+1}\delta^i_j\Big(R_{[mn]}-
    \big(\sigma_{[m\alpha|n]}+\sigma_{[m\beta}L^\beta_{\underset\vee{\alpha
    n}]}+\sigma_{[m\alpha}\underset2\nu{}_{n]}\big)\varphi^\alpha\Big)\\&
    +\sigma_{j[m|n]}\varphi^i+\sigma_{j[m}\sigma_{\alpha n]}\varphi^\alpha\varphi^i
    +\sigma_{j[m}L^i_{\underset\vee{\alpha n}]}\varphi^\alpha+
    \sigma_{j[m}\underset2\nu{}_{n]}\varphi^i-\delta^i_{[m}\underset2\mu\sigma_{jn]}\\&-
    \dfrac1{N+1}\Big(\delta^i_{[m}L^\alpha_{\underline{j\alpha}|n]}+\big(\delta^i_{[m}\sigma_{j\alpha|n]}+
    \delta^i_{[m}\sigma_{j\beta}L^\beta_{\underset\vee{\alpha
    n}]}+\delta^i_{[m}\sigma_{j\alpha}\underset2\nu{}_{n]}\big)\varphi^\alpha+
    \delta^i_{[m}\underset2\mu\sigma_{jn]}\Big)\\&+
    \dfrac1{N+1}\delta^i_{[m}\sigma_{jn]}\big(L^\beta_{\underline{\alpha\beta}}+\sigma_{\alpha\beta}\varphi^\beta\big)\varphi^\alpha
    \\&-\dfrac1{(N+1)^2}\big(L^\alpha_{\underline{j\alpha}}+\sigma_{j\alpha}\varphi^\alpha\big)
    \big(\delta^i_{[m}L^\beta_{\underline{n]\beta}}+\delta^i_{[m}\sigma_{n]\beta}\varphi^\beta\big),
    \endaligned\label{eq:Weylbasicsimpi32}\\\displaybreak[0]
    &\aligned
    \underset2{\widetilde
    W}{}^i_{jmn}&=R^i_{jmn}+\dfrac1{N+1}\delta^i_jR_{[mn]}+\dfrac
    N{N^2-1}\delta^i_{[m}R_{jn]}+\dfrac1{N^2-1}\delta^i_{[m}R_{n]j}
    \\&+\sigma_{j[m|n]}\varphi^i+\sigma_{j[m}L^i_{\underset\vee{\alpha
    n}]}\varphi^\alpha+\sigma_{j[m}\underset2\nu{}_{n]}\varphi^i+
    \sigma_{j[m}\sigma_{\alpha
    n]}\varphi^\alpha\varphi^i\\&+
    \dfrac1{N-1}\big(\delta^i_{[m}\sigma_{jn]|\alpha}+\delta^i_{[m}\sigma_{jn]}L^\beta_{\underset\vee{\alpha\beta}}+
    \delta^i_{[m}\sigma_{jn]}\underset2\nu{}_\alpha\big)\varphi^\alpha\\&+
    \dfrac1{N-1}\big(\delta^i_{[m}\sigma_{jn]}\sigma_{\alpha\beta}-\delta^i_{[m}\sigma_{j\alpha}\sigma_{n]\beta}\big)
    \varphi^\alpha\varphi^\beta\\&
    -\dfrac N{N^2-1}
    \big(\delta^i_{[m}\sigma_{\alpha j|n]}+\delta^i_{[m}\sigma_{\beta j}L^\beta_{\underset\vee{\alpha
    n}]}+\delta^i_{[m}\sigma_{\alpha j}\underset2\nu{}_{n]}\big)\varphi^\alpha\\&
    -\dfrac1{N^2-1}
    \big(\delta^i_{[m}\sigma_{\alpha n]|j}+\delta^i_{[m}\sigma_{\beta n]}L^\beta_{\underset\vee{\alpha
    j}}+\delta^i_{[m}\sigma_{\alpha
    n]}\underset2\nu{}_{j}\big)\varphi^\alpha,
    \endaligned\label{eq:basicWeylpi32}\\\displaybreak[0]
    &\aligned
    \underset2{\mathcal W}{}^i_{0.(p^1).(p^2).jmn}&=\underset2{\widetilde{\mathcal W}}{}^i_{jmn}
    +uL^i_{\underset\vee{jm}|n}+u'L^i_{\underset\vee{jn}|m}\\&-u\big(\omega^i_{(p_1^1).\alpha n}L^\alpha_{\underset\vee{jm}}-
    \omega^\alpha_{(p_2^1).jn}L^i_{\underset\vee{\alpha m}}-
    \omega^\alpha_{(p_3^1).mn}L^i_{\underset\vee{j\alpha}}\big)\\&-u'\big(
    \omega^i_{(p_1^2).\alpha m}L^\alpha_{\underset\vee{jn}}-
    \omega^\alpha_{(p_2^2).jm}L^i_{\underset\vee{\alpha n}}-
    \omega^\alpha_{(p_3^2).mn}L^i_{\underset\vee{j\alpha}}\big),
    \endaligned\label{eq:inv1KtoKpi32}\\\displaybreak[0]
    &\aligned
    \underset2W{}^i_{0.(p^1).(p^2).jmn}&=\underset2{\widetilde{W}}{}^i_{jmn}
    +uL^i_{\underset\vee{jm}|n}+u'L^i_{\underset\vee{jn}|m}\\&-u\big(\omega^i_{(p_1^1).\alpha n}L^\alpha_{\underset\vee{jm}}-
    \omega^\alpha_{(p_2^1).jn}L^i_{\underset\vee{\alpha m}}-
    \omega^\alpha_{(p_3^1).mn}L^i_{\underset\vee{j\alpha}}\big)\\&-u'\big(
    \omega^i_{(p_1^2).\alpha m}L^\alpha_{\underset\vee{jn}}-
    \omega^\alpha_{(p_2^2).jm}L^i_{\underset\vee{\alpha n}}-
    \omega^\alpha_{(p_3^2).mn}L^i_{\underset\vee{j\alpha}}\big),
    \endaligned\label{eq:inv2KtoKpi32}
  \end{align}

  \noindent for $p^1_1,\ldots,p^2_3\in\{1,2\}$,
  $\omega^i_{(1).jk}=L^i_{\underline{jk}}$ and
  $\omega^i_{(2).jk}=\omega^i_{jk}$, for the geometrical object $\omega^i_{jk}$
  given in the equation (\ref{eq:pi31omega}).

  The next theorem holds.

  \begin{thm}
    Let $f:\mathbb{GA}_N\to\mathbb{G\overline A}{}_N$ be an equitorsion almost
    geodesic mapping of the type $\underset2\pi{}_3$.

    The geometrical object $\underset2{\widetilde{\mathcal T}}{}^i_{jk}$ given
    by \emph{(\ref{eq:basicThomaspi32})} is the basic invariant of
    the Thomas type for the mapping $f$. The invariance of this
    geometrical object is total.

    The geometrical object $\underset2{\widetilde{\mathcal W}}{}^i_{jmn}$ given
    by \emph{(\ref{eq:Weylbasicsimpi32})} is the basic associated
    invariant for the mapping $f$. The invariance of this
    geometrical object is valued. It is total if and only if the
    mapping $f$ has the property of reciprocity.

    The geometrical object $\underset2{\widetilde W}{}^i_{jmn}$ given by
    \emph{(\ref{eq:basicWeylpi32})} is the associated derived invariant of the
    Weyl type for the mapping $f$.  The invariance of this
    geometrical object is valued. It is total if and only if the
    mapping $f$ has the property of reciprocity.

    The geometrical objects $\underset2{\mathcal W}{}^i_{jmn}$,
    $\underset2W{}^i_{jmn}$, given by
    \emph{(\ref{eq:inv1KtoKpi32}, \ref{eq:inv2KtoKpi32})}, are the invariants
    for the equitorsion third type almost geodesic mapping $f$.  The invariance of
    these
    geometrical objects are valued. They are total if and only if the
    mapping $f$ has the property of reciprocity.
    \qed
  \end{thm}

  \section{Conclusion}

  We obtained novel invariants for the almost geodesic mappings of
  the third type of a non-symmetric affine connection space in this
  paper.

  In the Section 3, the invariants for equitorsion almost geodesic
  mappings are presented. The method used for obtaining of these
  invariants \big(see \cite{jainv1}\big), simplified the
  corresponding method presented by Sinyukov \cite{sinjukov} and
  used latter in \cite{jaljubicamicapi3, micapi3spec}.


  The results obtained in this paper motivate the authors to continue their research
  about invariants for almost geodesic mappings of
  non-symmetric affine connection spaces.

  \section*{Acknowledgements}

  The authors thank to the anonymous referee who estimated the
  quality of this paper.


\begin{thebibliography}{33}

    \bibitem{mikpi1232014} \textbf{V. Berezovski, J. Mike\v s},
    \emph{Almost Geodesic Mappings of Affine Connection Spaces},
    (in Russian), Itogi Nauki i Tehniki, 126 (2014), 62--95.

    \bibitem{mikespi1232015} \textbf{V. Berezovski, J. Mike\v s},
    \emph{Almost Geodesic Mappings of Spaces With Affine
    Connection}, Journal of Mathematical Sciences, Vol. 207 (2015), No. 3,
    389--409.

    \bibitem{mikespi22020} \textbf{V. Berezovski, J. Mike\v s, L.
    R\'yparov\'a, A. Sabykanov}, \emph{On Canonical Almost Geodesic
    Mappings of Type $\pi_2(e)$}, Mathematics 2020, 8, 54;
    doi:10.3390/math8010054.


    \bibitem{eisNRG} \textbf{L. P. Eisenhart}, \emph{Non-Riemannian
Geometry}, New York, 1927.

    \bibitem{miknovi2019} \textbf{J. Mike\v s, S. B\'asc\'o, et
al.}, \emph{Differential geometry of special mappings}, Palacky
Univ., 2019.

\bibitem{mincic4} \textbf{S. M. Min\v ci\'c},
\emph{Curvature tensors of the space of non-symmetric affine
connexion, obtained from the curvature pseudotensors}, Matemati\v
cki Vesnik, 13 (28) (1976), 421--435.

\bibitem{mincic2} \textbf{S. M. Min\v ci\'c}, \emph{Independent
curvature tensors and pseudotensors of spaces with non-symmetric
affine connexion}, Coll. Math. Soc. J\'anos Bolayai, 31. Dif. geom.,
Budapest (Hungary), (1979), 445--460.

\bibitem{mincic1} \textbf{S. M. Min\v ci\'c}, \emph{On Ricci
Type Identities in Manifolds With Non-Symmetric Affine Connection},
PUBLICATIONS DE L'INSTITUT MATH\'EMATIQUE Nouvelle s\'erie, tome 94
(108) (2013), 205--217.

\bibitem{physinvss} \textbf{\DJ. Mu\v sicki, B. Mili\'c},
\emph{Mathematical Foundations of Theoretical Physics With a
Collection of Solved Problems}, (in Serbian), University of
Belgrade, 1974.

    \bibitem{petrovicmiskolc} \textbf{M. Z. Petrovi\'c},
    \emph{Canonical Almost Geodesic Mappings of Type
    $\underset\theta\pi{}_2(0,F)$, $\theta\in\{1,2\}$, Between
    Generalized Parabolic K\"ahler Manifolds}, Miskolc Mathematical Notes,
    Vol. 19 (2018), No. 1, pp. 469--482.

    \bibitem{petrovicmicapi1} \textbf{M. Z. Petrovi\'c, M. S.
    Stankovi\'c}, \emph{Special Almost Geodesic Mappings of the First Type of Non-symmetric Affine Connection Spaces},
    Bull. Malays. Math. Sci. Soc. 40, 1353--1362 (2017).

    \bibitem{sinjukov} \textbf{N. S. Sinyukov}, \emph{Geodesic mappings
    of Riemannian spaces},  (in Russian), "Nauka", Moscow, 1979.
    \bibitem{mikpi21999} \textbf{V. S. Sobchuk, J. Mike\v s, O.
    Pokorn\'a}, \emph{On Almost Geodesic Mappings $\pi_2$ Between
    Semisymmetric Riemannian Spaces}, Novi Sad J. Math., Vol. 29,
    No. 3, 1999, 309--312.

    \bibitem{micapi1} \textbf{M. S. Stankovi\'c}, \emph{First Type
    Almost Geodesic Mappings of General Affine Connection Spaces},
    Novi Sad J. Math., Vol. 29, No. 3, 1999, 313--323.

    \bibitem{micapi2} \textbf{M. S. Stankovi\'c}, \emph{On a Canonic
    Almost Geodesic Mappings of the Second Type of Affine Spaces},
    Filomat, 13 (1999), 105--114.

    \bibitem{micapi3} \textbf{M. S. Stankovi\'c}, \emph{On a Special
    Almost Geodesic Mappings of Third Type of Affine Spaces}, Novi
    Sad J. Math., Vol. 31, No. 2, 2001, 125--135.

    \bibitem{micapi3spec} \textbf{M. S. Stankovi\'c}, \emph{Special
    equitorsion almost geodesic mappings of the third type of
    non-symmetric affine connection spaces}, Applied Mathematics and
    Computation 244 (2014), 695--701.

    \bibitem{jainv1} \textbf{N. O. Vesi\'c}, \emph{Basic Invariants of
    Geometric Mappings}, Miskolc Mathematical Notes, Vol. 21 (2020), No.
    1, pp. 473--487.

    \bibitem{jaljubicamicapi3} \textbf{N. O. Vesi\'c, Lj. S.
    Velimirovi\'c, M. S. Stankovi\'c}, \emph{Some Invariants of
    Equitorsion Third Type Almost Geodesic Mappings},
    Mediterr. J. Math. 13 (2016), 4581--4590.

  \end{thebibliography}
\end{document}